\documentclass[final,leqno,onefignum,onetabnum]{siamltex1213}

\usepackage{tikz}
\usetikzlibrary{arrows,chains,matrix,positioning,scopes}
\makeatletter
\tikzset{join/.code=\tikzset{after node path={%
\ifx\tikzchainprevious\pgfutil@empty\else(\tikzchainprevious)%
edge[every join]#1(\tikzchaincurrent)\fi}}}
\makeatother

\tikzset{>=stealth',every on chain/.append style={join},
         every join/.style={->}}
\tikzstyle{labeled}=[execute at begin node=$\scriptstyle,
   execute at end node=$]

\usepackage[hyphenbreaks]{breakurl}
\usepackage{cite}
\usepackage{amsmath, amssymb, amsfonts}
\usepackage{array}

\allowdisplaybreaks

\usepackage{rotating}
\usepackage{tablefootnote}

\usepackage{algorithm2e}
\usepackage[noend]{algpseudocode}

\usepackage[multiple]{footmisc}

\title{Counterexample to global convergence of DSOS and SDSOS hierarchies}

\author{C\'edric Josz\footnotemark[1] }

\usepackage{soul}

\begin{document}
\maketitle

\renewcommand{\thefootnote}{\fnsymbol{footnote}}

\footnotetext[1]{Laboratory for Analysis and Architecture of Systems (LAAS), French National Center for Scientific Research (CNRS), 7, avenue du Colonel Roche, Toulouse, 31000, France (\email{cedric.josz@gmail.com}). The research was funded by the European Research Council (ERC) under the European Union's Horizon 2020 research and innovation program (grant agreement 666981 TAMING).}

\renewcommand{\thefootnote}{\arabic{footnote}}

\slugger{mms}{xxxx}{xx}{x}{x--x}

\begin{abstract}
We exhibit a convex polynomial optimization problem for which the diagonally-dominant sum-of-squares (DSOS) and the scaled diagonally-dominant sum-of-squares (SDSOS) hierarchies, based on linear programming and second-order conic programming respectively, do not converge to the global infimum. The same goes for the r-DSOS and r-SDSOS hierarchies. This refutes the claim in the literature according to which the DSOS and SDSOS hierarchies can solve any polynomial optimization problem to arbitrary accuracy. In contrast, the Lasserre hierarchy based on semidefinite programming yields the global infimum and the global minimizer with the first order relaxation. We further observe that the dual to the SDSOS hierarchy is the moment hierarchy where every positive semidefinite constraint is relaxed to all necessary second-order conic constraints. As a result, the number of second-order conic constraints grows quadratically in function of the size of the positive semidefinite constraints in the Lasserre hierarchy. Together with the counterexample, this suggests that DSOS and SDSOS are not necessarily more tractable alternatives to sum-of-squares.
\end{abstract}



\pagestyle{myheadings}
\thispagestyle{plain}
\markboth{C\'EDRIC JOSZ}{Counterexample to DSOS and SDSOS hierarchies}
\section{Introduction}
\label{sec:Introduction}
Consider a polynomial optimization problem
\begin{equation}
\begin{array}{ll}
\inf_x & f(x) := \sum_\alpha f_{\alpha} x^\alpha \\
\mathrm{s.t.} & g_i(x) := \sum_\alpha g_{i,\alpha} x^\alpha \geqslant 0, \quad i=1,\ldots,m
\end{array}
\end{equation}
where we use the multi-index notation $x^\alpha := x_1^{\alpha_1} \cdots x_n^{\alpha_n}$ for $x \in {\mathbb R}^n$,
$\alpha \in {\mathbb N}^n$, and
where the data are polynomials $f, g_1, \ldots, g_m \in {\mathbb R}[x]$
so that in the above sums only a finite number of coefficients
$f_{\alpha}$ and $g_{i,\alpha}$ are nonzero. We will use the notation $|\alpha|:= \sum_{k=1}^n \alpha_k$.

In 2001, the Lasserre hierarchy \cite{lasserre-2000,lasserre-2001} (see also \cite{parrilo-2000b,parrilo-2003}) was proposed to find global solutions to polynomial optimization problems. It is also known as moment/sum-of-squares hierarchy in reference to the primal moment hierarchy and the dual sum-of-squares hierarchy. Its global convergence is guaranteed by Putinar's Positivstellensatz \cite{putinar-1993} proven in 1993. Typically, if one of the constraints is a ball $x_1^2 + \hdots + x_n^2 \leqslant 1$, then the sequence of lower bounds provided by the hierarchy converges to the global infimum of the polynomial optimization problem. In addition, there is zero duality at all relaxation orders \cite{josz-2015}.

The moment problem of order $d$ is defined as
\begin{equation}
\begin{array}{ll}
\inf_y & L_y(f)  \\
\mathrm{s.t.} & y_0 = 1 \\
 & M_d(y) \succcurlyeq 0 \\
 &  M_{d-k_i}(g_iy) \succcurlyeq 0, \quad i=1,\ldots,m
\end{array}
\end{equation}
where $\succcurlyeq 0$ denotes positive semidefiniteness, and where the Riesz functional, the moment matrix, and the localizing matrices are respectively defined by
\begin{equation}
\begin{array}{l}
L_y(f) := \sum_{\alpha} f_{\alpha} y_\alpha\\[.2cm]
M_d(y) := (y_{\alpha+\beta})_{|\alpha|,|\beta|\leqslant d} \\[.2cm]
M_{d-k_i}(g_iy) := (\sum_{\gamma} g_{i,\gamma} y_{\alpha+\beta+\gamma})_{|\alpha|,|\beta|\leqslant d-k_i}\\[.2cm]
k_i := \max\{ \lceil|\alpha|/ 2 \rceil ~\text{s.t.}~ g_{i,\alpha} \neq 0 \}.
\end{array}
\end{equation}
Above, $\lceil . \rceil$ denotes the ceiling of a real number.

The sum-of-squares problem of order $d$ is defined as
\begin{equation}
\begin{array}{ll}
\sup_{\lambda,\sigma} & \lambda  \\
\mathrm{s.t.} & f - \lambda = \sigma_0 + \sum_{k=1}^m \sigma_k g_k \\
 & \lambda \in \mathbb{R}, \sigma_0 \in \Sigma_d[x], \\
 & \sigma_i \in \Sigma_{d-k_i}[x], \quad i=1,\ldots,m.
\end{array}
\end{equation}
A polynomial $\sigma(x) = \sum_{|\alpha|\leqslant 2d} \sigma_{\alpha} x^{\alpha}$ is a sum-of-squares, i.e. it belongs to $\Sigma_d[x]$, if it is of the form
\begin{equation}
\sigma(x) = \sum_k \left( \sum_{|\alpha|\leqslant d} p_{k,\alpha} x^\alpha \right)^2 ~~~ \text{where}~~~ p_{k,\alpha} \in \mathbb{R}.
\end{equation}
This is equivalent to the existence of $(\varphi_{\alpha,\beta})_{|\alpha|,|\beta|\leqslant d} \succcurlyeq 0$ such that $\sum_{|\alpha|\leqslant 2d} \sigma_{\alpha} x^{\alpha} = \sum_{|\alpha|,|\beta|\leqslant d} \varphi_{\alpha,\beta} x^{\alpha+\beta}$.

Following \cite{ahmadi-cdc-2014,ahmadi-2017}, a polynomial $\sigma(x) = \sum_{|\alpha|\leqslant 2d} \sigma_{\alpha} x^{\alpha}$ is a scaled diagonally-dominant sum-of-squares (SDSOS) if it is of the form
\begin{equation}
\label{eq:sdsos}
\sigma(x) = \sum_k \left( p_k x^{\alpha(k)} + q_k x^{\beta(k)} \right)^2 ~~~ \text{where}~~~ \begin{array}{c} \scriptsize \alpha(k), \beta(k) \in \mathbb{N}^n, \\ p_k, q_k \in \mathbb{R}. \end{array}
\end{equation}
This is equivalent to the existence of $(\varphi_{\alpha,\beta})_{|\alpha|,|\beta|\leqslant d}$ such that $\sum_{|\alpha|\leqslant 2d} \sigma_{\alpha} x^{\alpha} = \sum_{|\alpha|,|\beta|\leqslant d} \varphi_{\alpha,\beta} x^{\alpha+\beta}$ where $\varphi$ is of the form
\begin{equation}
\label{eq:sum}
\sum_{\small \begin{array}{c} |\alpha|,|\beta| \leqslant d \\ \alpha \neq \beta \end{array}} \phi^{\alpha,\beta}
\end{equation}
where each matrix $\phi^{\alpha,\beta}$ has zeros everywhere expect for four entries $\phi^{\alpha,\beta}(\alpha,\alpha)$, $\phi^{\alpha,\beta}(\alpha,\beta)$, $\phi^{\alpha,\beta}(\beta,\alpha)$, $\phi^{\alpha,\beta}(\beta,\beta) \in \mathbb{R}$ such that
\begin{equation}
\label{eq:socp}
\begin{pmatrix}
\phi^{\alpha,\beta}(\alpha,\alpha)~  & ~\phi^{\alpha,\beta}(\alpha,\beta) \\
\phi^{\alpha,\beta}(\beta,\alpha)~ & ~\phi^{\alpha,\beta}(\beta,\beta)
\end{pmatrix}
\succcurlyeq 0.
\end{equation}
These can be viewed as second-order conic constraints. The equivalence between \eqref{eq:sdsos} and \eqref{eq:sum}-\eqref{eq:socp} is immediate and does not require \cite[Theorem 7,Theorem 8, Lemma 9]{ahmadi-2017}.

In the sum-of-squares problem of order $d$, if we restrain the sum-of-squares variables $\sigma_0, \hdots , \sigma_m$ to be SDSOS, then in the moment problem of order $d$, we are relaxing each semidefinite constraint as follows
\begin{equation}
\begin{pmatrix}
M_d(y)(\alpha,\alpha)~  & ~M_d(y)(\alpha,\beta) \\
M_d(y)(\beta,\alpha)~ & ~ M_d(y)(\beta,\beta)
\end{pmatrix}
\succcurlyeq 0 ~,~~~ \forall |\alpha|,|\beta| \leqslant d, ~~ \alpha \neq \beta,
\end{equation}
\begin{equation}
\begin{pmatrix}
M_{d-k_i}(g_iy)(\alpha,\alpha)~  & ~M_{d-k_i}(g_iy)(\alpha,\beta) \\
M_{d-k_i}(g_iy)(\beta,\alpha)~ & ~ M_{d-k_i}(g_iy)(\beta,\beta)
\end{pmatrix}
\succcurlyeq 0 ~,~~~ \forall |\alpha|,|\beta| \leqslant d-k_i, ~~ \alpha \neq \beta.
\end{equation}
Naturally, these can also be viewed as second-order conic constraints. The number of second-order conic constraints is \footnote{We use the notation
\begin{equation}
\begin{pmatrix}
n \\
k
\end{pmatrix} := \frac{n(n-1)\hdots(n-k+1)}{k!} 
\end{equation}
for integers $1 \leqslant k \leqslant n$ and where $k!$ stands for factorial.}
\begin{equation}
\begin{pmatrix}
\begin{pmatrix}
n+d\\ 
d
\end{pmatrix}\\
2
\end{pmatrix}
~~+~~
\sum_{i=1}^m
\begin{pmatrix}
\begin{pmatrix}
n+d-k_i\\ 
d-k_i
\end{pmatrix}\\
2
\end{pmatrix}
~~=~~ O(n^{2d})
\end{equation}
as the number of variables $n$ grows, to be compared with the semidefinite constraints of size up to $O(n^d)$ in the Lasserre hierarchy. In particular, with $n=10$ variables at order $d=2$, SDSOS relaxes the moment matrix of size $66 \times 66$ to $2,145$ second-order conic constraints. 

In the sum-of-squares problem of order $d$, if we restrain the sum-of-squares variables $\sigma_0, \hdots , \sigma_m$ to be DSOS, then in the moment problem of order $d$, we are relaxing each semidefinite constraint to the following linear constraints:\\\\
$\forall |\alpha|,|\beta| \leqslant d, ~~ \alpha \neq \beta,$
\begin{equation}
\begin{array}{c}
M_d(y)(\alpha,\alpha)\geqslant 0,~ M_d(y)(\beta,\beta) \geqslant 0, \\ 
M_d(y)(\alpha,\alpha) + 2M_d(y)(\alpha,\beta) + M_d(y)(\beta,\beta) \geqslant 0, \\
M_d(y)(\alpha,\alpha) - 2M_d(y)(\alpha,\beta) + M_d(y)(\beta,\beta) \geqslant 0,
\end{array}
\end{equation}
$\forall |\alpha|,|\beta| \leqslant d-k_i, ~~ \alpha \neq \beta,$
\begin{equation}
\begin{array}{c}
M_{d-k_i}(g_iy)(\alpha,\alpha)\geqslant 0,~ M_{d-k_i}(g_iy)(\beta,\beta) \geqslant 0, \\ 
M_{d-k_i}(g_iy)(\alpha,\alpha) + 2M_{d-k_i}(g_iy)(\alpha,\beta) + M_{d-k_i}(g_iy)(\beta,\beta) \geqslant 0, \\
M_{d-k_i}(g_iy)(\alpha,\alpha) - 2M_{d-k_i}(g_iy)(\alpha,\beta) + M_{d-k_i}(g_iy)(\beta,\beta) \geqslant 0.
\end{array}
\end{equation}
In the language of the lift-and-project method \cite{laurent-2003}, the three above equations correspond to the lifting of $(x^\alpha)^2g_i(x) \geqslant 0, (x^\beta)^2g_i(x) \geqslant 0, (x^\alpha + x^\beta)^2 g_i(x) \geqslant 0$ and $(x^\alpha - x^\beta)^2g_i(x) \geqslant 0$. For a nice description of the DSOS and SDSOS hierarchies, see  \cite[(QM-DSOSr),(QM-SDSOSr)]{kuang-2017} and \cite{kuang-2017bis}.\\

We now turn our attention to claims that have been made in the literature:\\
\begin{enumerate}
\item In \cite{ahmadi-2014}, it is stated that:\\\\
\textit{``A particularly nice feature of DSOS and SDSOS opti­mization is that they enjoy many of the same theoretical guarantees that underly SOS optimization.''}\\\\
 ``\textit{Theorem 2.5 (Ahmadi, Majumdar, 13): Consider the general polynomial optimization problem (POP) in 1. There is a hierarchy of linear programs based on optimization over dsos polynomials that can solve POP to arbitrary accuracy.}

\textit{
This theorem is similar to the Parrilo and Lasserre hierarchies for polynomial optimization that are instead based on a search over sum of squares polynomials via semidefinite programming [11], [13].}''\\

\item In \cite{ahmadi-2013}, it is stated that:\\\\
``\textit{Similar to the Lasserre/Parrilo hierarchies of SDP (based on Putinar, Schm\"udgen or Stengle's Positivstellensatz) that solve POP to global optimality,
our converse results imply that POP can be solved to global optimality using hierarchies of LP and SOCP coming from dsos and sdsos.}''\\
\end{enumerate}

The object of this manuscript is to propose a counterexample to these claims. To the best of our understanding, the response \cite{response} to this manuscript does not show how the converse results presented in \cite{ahmadi-2014,ahmadi-2013} imply that the counterexample can be solved to global optimality using hierarchies of linear programs and second-order conic programs.

The manuscript is organized as follows. Section \ref{sec:Counterexample} states the counterexample and lists its main properties. Section \ref{sec:Proof} contains the proof that the example is indeed a counterexample. Section \ref{sec:Discussion} discusses the consequences of our finding. Section \ref{sec:Conclusion} concludes the manuscript.

\section{Counterexample}
\label{sec:Counterexample}
Consider the following optimization problem:\\\\
\begin{equation}
\boxed{ \inf_{x_1,x_2 \in \mathbb{R}} ~~~ (-2+x_1+x_2)^2 ~~~~~\text{subject to}~~~~~ 1-x_1^2-x_2^2 \geqslant 0 }
\end{equation}
\\\\
This problem is endowed with the following properties:\\
\begin{enumerate}
\item it is a convex optimization problem;\\
\item the objective and constraint functions are sum-of-squares-convex;\\
\item the feasible set is compact;\\
\item Slater's condition is satisfied;\\
\item the quadratic module associated to the constraint is Archimedean.\\\\
\end{enumerate}
Based on the above properties, the Karush Kuhn Tucker conditions hold at optimality
\begin{equation}
\left\{
\begin{array}{c}
-4+2x_1+2x_2 + 2\lambda x_1 = 0 \\[.1cm]
-4+2x_1+2x_2 + 2\lambda x_2 = 0 \\[.1cm]
x_1^2 + x_2^2 - 1 \leqslant 0, ~ \lambda \geqslant 0 \\[.1cm]
\lambda (x_1^2 + x_2^2 - 1) = 0
\end{array}
\right.
\end{equation}
yielding the unique primal-dual global solution
\begin{equation}
x_1=x_2 = \frac{1}{\sqrt{2}} ~~~\text{and}~~~ \lambda = 2(\sqrt{2}-1)
\end{equation}
and the global value 
\begin{equation}
2(3-2\sqrt{2}).
\end{equation}
\\
We will see in the next section that the Lasserre hierarchy finds the global minimizer and the global infimum at the first order relaxation. It also finds the Karush Kuhn Tucker multiplier. In contrast, the SDSOS hierarchy does not globally converge. In fact, it has a relaxation gap equal to 2 at all orders. As a consequence, the weaker DSOS hierarchy has a relaxation gap of at least 2 at all orders, and does not globally converge either.

\section{Proof}
\label{sec:Proof}
For convenience, let $ f(x_1,x_2) : = (-2+x_1+x_2)^2 $ denote the objective function and let $g(x_1,x_2) := 1-x_1^2-x_2^2$ denote the constraint function. The Hessians satisfy
\begin{equation}
\nabla^2 f(x_1,x_2) = \begin{pmatrix} \hphantom{-}2 & -2 \\ -2 & \hphantom{-}2 \end{pmatrix} \succcurlyeq 0,
\end{equation}
\begin{equation}
-\nabla^2 g(x_1,x_2) = \begin{pmatrix} 1 & 0 \\ 0 & 1 \end{pmatrix} \succcurlyeq 0. ~~~~~~
\end{equation}
According to \cite[Definition 2.3]{lasserre-2010}, $f$ and $g$ are sum-of-squares-convex. Applying \cite[Theorem 5.15]{lasserre-2010}, we deduce that the Lasserre hierarchy globally converges at the first order relaxation. Let's double check: a primal-dual feasible point to the first order relaxation of the Lasserre hierarchy is given by
\begin{equation}
M_1(y) = 
\begin{array}{cccc}
 & 1 & x_1 & x_2 \\[.25em]
1 & 1/\sqrt{1} & \hphantom{-}1/\sqrt{2} & \hphantom{-}1/\sqrt{2} \\[.25em]
x_1 & 1/\sqrt{2} & \hphantom{-}1/\sqrt{4} & \hphantom{-}1/\sqrt{4} \\[.25em]
x_2 & 1/\sqrt{2} & \hphantom{-}1/\sqrt{4} & \hphantom{-}1/\sqrt{4} \\[.25em]
\end{array}
\end{equation}
and 
$$ f(x_1,x_2) - 2(3-2\sqrt{2}) = (\sqrt{2}-1) (x_1 - x_2)^2 + \sqrt{2} (-\sqrt{2} + x_1 + x_2)^2 + 2(\sqrt{2}-1)g(x_1,x_2). $$
The values of the above primal and the dual points are both equal to $2(3-2\sqrt{2})$, making it the primal-dual optimal value of the first order relaxation. Moreover, the rank of the moment matrix is equal to one, hence $2(3-2\sqrt{2})$ is the global value of the polynomial problem and $x_1=x_2 = \frac{1}{\sqrt{2}}$ is a global solution. According to \cite[Theorem 5.12 (c)]{lasserre-2010}, the Karush Kuhn Tucker multiplier associated to the constraint is given by evaluating at the global solution the sum-of-squares associated to the constraint, in this case the constant polynomial $2(\sqrt{2}-1)$.\\\\
A primal-dual feasible point to the first order relaxation of the SDSOS hierarchy is given by
\begin{equation}
M_1(y) = 
\begin{array}{cccc}
 & 1 & x_1 & x_2 \\[.25em]
1 & 1/\sqrt{1} & \hphantom{-}1/\sqrt{2} & \hphantom{-}1/\sqrt{2} \\[.25em]
x_1 & 1/\sqrt{2} & \hphantom{-}1/\sqrt{4} & -1/\sqrt{4} \\[.25em]
x_2 & 1/\sqrt{2} & -1/\sqrt{4} & \hphantom{-}1/\sqrt{4} \\[.25em]
\end{array}
\end{equation}
and
$$ f(x_1,x_2) - 4(1-\sqrt{2}) =  \frac{\sqrt{2}}{2}( -\sqrt{2} + 2 x_1)^2 +  \frac{\sqrt{2}}{2}( -\sqrt{2} + 2 x_2)^2 + (x_1+ x_2)^2 + 2\sqrt{2}g(x_1,x_2).$$
The values of the above primal and the dual points are both equal to $4(1-\sqrt{2})$, making it the primal-dual optimal value of the first order relaxation. We now show that it is the primal-dual optimal value of the SDSOS relaxation of any order. The sequence defined by
\begin{equation}
\label{eq:sequence}
y_{(\alpha_1,\alpha_2)} := \frac{ 1 + (-1)^{\alpha_1} + (-1)^{\alpha_2} -(-1)^{\alpha_1+\alpha_2}}{2\sqrt{2^{\alpha_1+\alpha_2}}} ~,~~~ \forall \alpha_1,\alpha_2 \in  \mathbb{N},
\end{equation}
that is to say
\begin{equation}
M(y) =
\begin{array}{cccccccc}
 & 1 & x_1 & x_2 & x_1^2 & x_1 x_2 & x_2^2 & \\[.25em]
1 & \hphantom{-}1/\sqrt{1} & \hphantom{-}1/\sqrt{2} & \hphantom{-}1/\sqrt{2} & \hphantom{-}1/\sqrt{4} & -1/\sqrt{4} & \hphantom{-}1/\sqrt{4} & \hdots \\[.25em]
x_1 & \hphantom{-}1/\sqrt{2} & \hphantom{-}1/\sqrt{4} & -1/\sqrt{4} & \hphantom{-}1/\sqrt{8} & \hphantom{-}1/\sqrt{8} & \hphantom{-}1/\sqrt{8} & \\[.25em]
x_2 & \hphantom{-}1/\sqrt{2} & -1/\sqrt{4} & \hphantom{-}1/\sqrt{4} & \hphantom{-}1/\sqrt{8} & \hphantom{-}1/\sqrt{8} & \hphantom{-}1/\sqrt{8} & \\[.25em]
x_1^2 & \hphantom{-}1/\sqrt{4} & \hphantom{-}1/\sqrt{8} & \hphantom{-}1/\sqrt{8} & \hphantom{-}1/\sqrt{16} & -1/\sqrt{16} & \hphantom{-}1/\sqrt{16} & \\[.25em]
x_1x_2 & -1/\sqrt{4} & \hphantom{-}1/\sqrt{8} & \hphantom{-}1/\sqrt{8} & -1/\sqrt{16} & \hphantom{-}1/\sqrt{16} & -1/\sqrt{16} & \\[.25em]
x_2^2 & \hphantom{-}1/\sqrt{4} & \hphantom{-}1/\sqrt{8} & \hphantom{-}1/\sqrt{8} & \hphantom{-}1/\sqrt{16} & -1/\sqrt{16} & \hphantom{-}1/\sqrt{16} & \\[.25em]
 & \vdots & & & & & & \ddots
\end{array}
\end{equation}
yields a primal optimal solution (after truncation) to the relaxation of any order of the SDSOS hierarchy because
\begin{enumerate}
\item $L_y(f) = 4(1-\sqrt{2});$\\
\item
$
\begin{pmatrix}
y_{2\alpha} & y_{\alpha+\beta} \\
y_{\alpha+\beta} & y_{2\beta}
\end{pmatrix} \succcurlyeq 0 ~,~~~ \forall \alpha,\beta \in  \mathbb{N}^2;\\
$
\item $y_{(\alpha_1,\alpha_2)} = y_{(\alpha_1+2,\alpha_2)} + y_{(\alpha_1,\alpha_2+2)}  ~,~~~ \forall \alpha_1,\alpha_2 \in  \mathbb{N}.$
\end{enumerate}
The last two points follow readily from the equation
\begin{equation}
|y_{\alpha}| = \frac{1}{\sqrt{2^{|\alpha|}}} ~,~~~ \forall \alpha \in \mathbb{N}^2
\end{equation}
and from the fact that $y_{\alpha}$ has a minus sign if and only if $\alpha_1$ and $\alpha_2$ are odd. 
This terminates the proof.
\section{Discussion}
\label{sec:Discussion}

First, we discuss the result which guarantees the convergence of the Lasserre hierarchy. 
\begin{theorem}[Putinar's Positivstellensatz \cite{putinar-1993}]
Assume that there exists $R>0$ and sums-of-squares $p_0,\hdots,p_m$ such that $R^2 - x_1^2 - \hdots - x_n^2 = p_0 + \sum_{i=1}^m p_i g_i $.
If $f>0$ on $\{ x \in \mathbb{R}^n ~|~ g_1(x) \geqslant 0 ~, \hdots~ ,~ g_m(x) \geqslant 0 \}$, then there exists sum-of-squares $\sigma_0,\hdots,\sigma_m$ such that
\begin{equation}
f = \sigma_0 + \sum_{i=1}^m \sigma_i g_i.
\end{equation}
\end{theorem}
In fact, the convergence of the Lasserre hierarchy is equivalent to Putinar's Positivstellensatz. If the SDSOS hierarchy were globally convergent as claimed in the literature, then under the same assumptions, the conclusion of Putinar's Positivstellensatz could be strenghened to the existence of SDSOS polynomials $\sigma_0,\hdots,\sigma_m$. In other words, $\sigma_0,\hdots,\sigma_m$ could be chosen to be sum-of-squares of polynomials of two terms at most.
As a byproduct of the previous section, there does not exist SDSOS polynomials $\sigma_0$ and $\sigma_1$ such that 
\begin{equation}
(-2+x_1+x_2)^2 = \sigma_0(x_1,x_2) + \sigma_1(x_1,x_2) ( 1 - x_1^2 - x_2^2 ) 
\end{equation}
even though $(-2+x_1+x_2)^2 > 0$ whenever $1 - x_1^2 - x_2^2 \geqslant 0$. As a consequence, the existence of the more stringent DSOS polynomials $\sigma_0$ and $\sigma_1$ is also impossible. In light of the simplicity of the counterexample, it is of little interest to investigate conditions under which one may seek SDSOS or DSOS polynomials in Putinar's Positivstellensatz.

In \cite[Definition 11]{ahmadi-2017}, the notions of r-DSOS and r-SDSOS are proposed to strenghen the DSOS and SDSOS hierarchies. They consist in premutliplying by $(x_1^2 + \hdots + x_n^2)^{r}$ where $r \in \mathbb{N}$ in the DSOS/SDSOS dual optimization problem. On our example this reads:
\begin{equation}
\label{eq:con}
\sup_{\lambda,\sigma} ~~\lambda ~~~~\text{s.t.}~~~~ (x_1^2+x_2^2)^{r} ( f(x_1,x_2) - \lambda ) =  \sigma_0(x_1,x_2) + \sigma_1(x_1,x_2) ( 1 - x_1^2 - x_2^2 ).
\end{equation}
We now show that the dual optimal value of the r-SDSOS hierarchy cannot exceed 0, which is strictly inferior to the global infimum $2(3-2\sqrt{2})$. As a consequence, the r-SDSOS and the r-DSOS hierarchies do not convergence globally.

The sequence $y$ defined in \eqref{eq:sequence} satisfies
\begin{equation}
\begin{array}{rcl}
L_y((x_1^2+x_2^2)^{r}f) & = & \sum\limits_{k=0}^r \begin{pmatrix} n \\ k \end{pmatrix} L_y( x_1^{2k}x_2^{2(n-k)}f ) \\[.4cm]
 & = & \sum\limits_{k=0}^r \begin{pmatrix} n \\ k \end{pmatrix} L_y( f ) / \sqrt{2^{2k+2(n-k)}} \\[.4cm]
 & = & L_y( f ) / 2^{n} \sum\limits_{k=0}^r \begin{pmatrix} n \\ k \end{pmatrix}  \\[.4cm]
 & = & L_y( f ).
\end{array}
\end{equation}
This means that the sequence $y$ remains primal optimal for all orders of the SDSOS hierarchy if we replace the objective function $f$ by $(x_1^2+x_2^2)^{r}f$. As a result, there does not exist SDSOS polynomials $\sigma_0$ and $\sigma_1$ such that 
\begin{equation}
(x_1^2+x_2^2)^{r}(-2+x_1+x_2)^2 = \sigma_0(x_1,x_2) + \sigma_1(x_1,x_2) ( 1 - x_1^2 - x_2^2 ) 
\end{equation}  
for all $r \in \mathbb{N}$. Hence, if $\lambda \geqslant 0$ in \eqref{eq:con}, then we have a contradiction.\\

Second, we discuss two claims that have been made in the literature regarding the tractability of DSOS and SDSOS hiearchies compared to the Lasserre hierarchy:\\
\begin{enumerate}
\item In the conclusion of \cite{ahmadi-2017}, it is stated that \\\\
``\textit{Our numerical examples from a diverse range of applications including polynomial optimization, combinatorial optimization, statistics and machine learning, derivative pricing, control theory, and robotics demonstrate that with reasonable tradeoffs in optimality, we can handle problem sizes that are well beyond the current capabilities of SOS programming. In particular, we have shown that our approach is able to tackle dense problems with as many as 70 polynomial variables (with degree-4 polynomials).}''\\
\item On page 3 of \cite{ahmadi-2017}, it is stated that\\\\
``\textit{POP contains as special cases many important problems in operations research; e.g., the optimal power flow problem in power engineering [39]''} ([39] = M. Huneault and F. Galiana, A survey of the optimal power flow literature, IEEE Transactions on Power Systems, 6 (1991), pp. 762-770.)\\
\end{enumerate}

The survey paper \cite{ahmadi-2017} titled \textit{``DSOS and SDSOS Optimization: More Tractable Alternatives to Sum of Squares and Semidefinite Optimization''} does not mention that in \cite{kuang-2017bis}, the DSOS and SDSOS hierarchies provide strict lower bounds to instances of the optimal power flow problem \cite{carpentier-1962} with up to 600 variables and 1,300 constraints. It also does not mention that in \cite{josz-molzahn-2018}, the \textit{multi-ordered} Lasserre hierarchy finds global minimizers to instances of the optimal power flow problem with up to 4,500 variables and 14,500 constraints. The instances correspond to the European high-voltage electricity transmission network. The data used was collected from 23 different national transmission system operators and is available in \cite{josz-2016}.\\



The multi-ordered Lasserre hierarchy \cite{josz-phd,cdc2015,molzahn-pscc2016,josz-molzahn-2018,mh_sparse_msdp} was proposed to exploit sparsity in general polynomial optimization problems. The main ideas are: 1) to use a different relaxation order for each constraint, and 2) to iteratively seek a closest measure to the truncated moment data until a measure matches the truncated data. The multi-ordered Lasserre builds on previous work on chordal sparsity \cite{waki-2006,gron1984,schweighofer-2007,kuhlmann-2007}. Proof of convergence of the multi-ordered hierarchy can be found in \cite[Section 6]{josz-molzahn-2018}. There exists other avenues of work to deal with large scale problems in polynomial optimization, such as the bounded sum-of-squares (BSOS) hierarchy \cite{toh-2017,weisser-2017} which is globally convergent. Also, the authors of \cite{riener-2013} paved the way for exploiting symmetry in the Lasserre hierarchy, which can reduce the computational burden.\\

Third and last, the idea of relaxing the semidefinite constraints in the Lasserre hierarchy to second-order conic constraints (as in the SDSOS hierarchy) was independently proposed in \cite{dan2015}. In that work, the moment constraint is maintained as a positive semidefinite constraint, but the localizing matrices are relaxed to multiple second-order conic constraints. This guarantees that the relaxation is stronger than the first order Lasserre relaxation. The paper considers only medium-sized instances of the optimal power flow problem (several hundreds of variables). In some instances, there is a computational gain, but in others, the approach is unable to find a minimizer whereas the Lasserre hierarchy does. In light of the counterexample of Section \ref{sec:Counterexample}, it is possible that the instances where the approach fails are also counterexamples.\\
 
\section{Conclusion}
\label{sec:Conclusion}
In conclusion, the DSOS, SDSOS, r-DSOS, and r-SDSOS hierarchies do not preserve the global convergence property associated with the Lasserre hierarchy. Rather, they constitute a heuristic based on the following idea: whenever one is faced with a positive semidefinite constraint in applied mathematics, one can always relax it to multiple linear constraints or multiple second-order conic constraints.

\bibliography{mybib}{}
\bibliographystyle{siam}

\end{document}